\documentclass{article}
\usepackage{amsmath}
\usepackage{amsthm}
\usepackage{amssymb}
\usepackage{graphicx}
\usepackage{hyperref}

\newtheorem{lemma}{Lemma}[section]
\newtheorem{proposition}[lemma]{Proposition}
\newtheorem{theorem}[lemma]{Theorem}

\newtheorem{definition}[lemma]{Definition}
\newtheorem{remark}[lemma]{Remark}

\def\row{\mathrm{row}}
\def\injection{\hookrightarrow}
\def\hom{\,\mathrm{Hom}}
\def\mod{\,\mathrm{mod}\,\,}
\def\aut{\mathrm{Aut}}
\def\semi{\rtimes}
\def\ints{\mathbb{Z}}

\def\sp{\mathrm{Sp}}

\def\mcg#1{\mathcal{M}_{#1}}
\def\handle#1{\mathcal{H}_{#1}}
\def\torelli#1{\mathcal{I}_{#1}}
\def\be{\begin{enumerate}}
\def\ee{\end{enumerate}}
\def\bi{\begin{itemize}}
\def\ei{\end{itemize}}

\title{Calculating the image of the second Johnson-Morita representation}

\author{
Joan S. Birman\footnote{The first author was supported in part by NSF grant  DMS-0405586.},
Tara E. Brendle\footnote{The second author was supported in part by NSF grant DMS-0606882.}
and Nathan Broaddus\footnote{The third author was supported in part by an NSF Postdoctoral Fellowship.}
}
\date{August 28, 2007}
\begin{document}

\maketitle

%Dedication ***************************
\begin{center}
\emph{\small This paper is dedicated, with respect, to Shigeyuki Morita.}
\end{center}
%end Dedication ***********************

\begin{abstract} Johnson has defined a surjective homomorphism from 
the {\em Torelli subgroup} of the mapping class group of the surface of 
genus $g$ with one boundary component to $\wedge^3 H$, the third 
exterior product of the homology of the surface.  Morita then
extended Johnson's homomorphism to a homomorphism from the entire 
mapping class group to $\frac{1}{2} \wedge^3 H \semi \sp(H)$.  
This {\em Johnson-Morita homomorphism} is not surjective, but its image is finite 
index in $\frac{1}{2} \wedge^3 H \semi \sp(H)$ \cite{Morita1993a}.  Here we give a 
description of the exact image of Morita's homomorphism.  
Further, we compute the image of the {\em handlebody subgroup} of the mapping class group under the same map.
\end{abstract}

\section{Introduction}
\begin{figure}[ht]
\begin{center}
\includegraphics[scale=0.32]{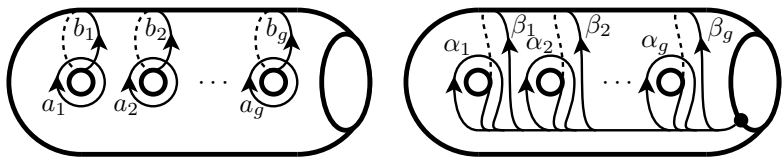}
\put(-280.0,17.4){$a_1$}
\put(-268.0,45.6){$b_1$}
\put(-252.8,17.4){$a_2$}
\put(-240.8,45.6){$b_2$}
\put(-220.8,25.6){$\cdots$}
\put(-205.0,17.2){$a_g$}
\put(-195.2,45.6){$b_g$}
\put(-128.0,40.0){$\alpha_1$}
\put(-102.4,47.0){$\beta_1$}
\put(-96.8,40.8){$\alpha_2$}
\put(-74.4,47.0){$\beta_2$}
\put(-68.0,25.6){$\cdots$}
\put(-53.2,41.6){$\alpha_g$}
\put(-28.0,46.0){$\beta_g$}
\caption{(a) A basis for $H_1(S_{g,1})$  \hspace{0.05in} (b) Generators for $\pi_1(S_{g,1})$\label{figure:curves}}
\end{center}
\end{figure}

Let $S_g$ be a closed surface of genus $g$.  We fix a 
closed disk $D$ in $S_{g}$, and by deleting its interior, obtain 
$S_{g,1}$, a genus $g$ surface with one boundary component, as 
illustrated in Figure~\ref{figure:curves}.  Let $\mcg{g}$ (resp. 
$\mcg{g,1}$) denote the mapping class group of the surface $S_g$ 
(resp. $S_{g,1}$).  In the case of $\mcg{g,1}$ we assume the 
boundary component is fixed pointwise.   

We choose a base point on $\partial S_{g,1}$, and let 
$\alpha_1, \ldots, \alpha_g,\beta_1, \ldots, \beta_g$ denote the 
based loops illustrated in Figure~\ref{figure:curves}(b).  Let 
$a_1,\dots, a_g,b_1,\dots, b_g$ denote the corresponding homology 
classes, as in Figure~\ref{figure:curves}(a).  It will sometimes be 
convenient to denote these same homology classes by 
$x_1, \dots, x_{2g}$ with the understanding that $x_i = a_i$ and 
$x_{i+g} = b_i$ for $1\leq i \leq g$.  Likewise, we will sometimes refer to the based loops $\alpha_1, \ldots, \alpha_g, \beta_1, \ldots, \beta_g$ by $\xi_1, \dots, \xi_{2g}$ with the understanding that $\xi_i = \alpha_i$ and 
$\xi_{i+g} = \beta_i$ for $1\leq i \leq g$.
 
Now, let $H = H_1(S_{g,1})$ be the free abelian group with generating set $\{a_1,\dots, a_g,b_1,\dots, b_g\}$ and $\pi = \pi_1(S_{g,1})$ which is a free group on the generating set $\{\alpha_1, \ldots, \alpha_g, \beta_1, \ldots, \beta_g\}$.  The action of $\mcg{g,1}$ on $\pi$ 
gives an injection $\mcg{g,1} \injection \aut(\pi)$.  More 
generally, we can compose with the homomorphism 
$\aut(\pi) \to \aut(\pi / \chi)$ for any characteristic subgroup 
$\chi \subset \pi$.  The \emph{lower central series} of the free 
group $\pi$ is a sequence of characteristic subgroups
defined inductively by setting $\pi^{(0)} = \pi$ and $\pi^{(k+1)} = 
[\pi, \pi^{(k)}]$.  We define the $k^{th}$ {\em Johnson-Morita 
representation} to be the map
\begin{equation*}
\rho_k : \mcg{g,1} \to \aut(\pi / \pi^{(k)})
\end{equation*}
We note that these maps were first studied by Johnson in 
\cite{Johnson1983b,Johnson1985} and subsequently developed by Morita 
in a series of papers \cite{Morita1993a,Morita1993b,Morita1996,Morita2001}.

Observe that the first Johnson-Morita map is just the classical 
symplectic representation 
$\rho_1 : \mcg{g,1} \to \sp(H)$ which is surjective (\cite{Burkhardt1890}, in particular pp. 209-212).  In \cite[Theorem 4.8]{Morita1993a} Morita 
shows that the image of $\rho_2$ is isomorphic to a subgroup of finite index 
in $\frac{1}{2} \wedge^3 H \semi \sp(H)$. 
Our first main result in this paper, 
given in Theorem~\ref{theorem:morita image}, is to identify the 
precise image $\rho_2(\mcg{g,1})$ using a formulation due to Perron 
\cite{Perron2004}. 

Let us now consider $S_g$ 
as $\partial X_g$, where $X_g$ is a genus $g$ handlebody.   Let $\handle{g}$ denote the {\em handlebody subgroup} of $\mcg{g}$, that is, the subgroup consisting of maps of $S_g$ which extend to the handlebody $X_g$.  There is a natural surjection $\mcg{g,1} \to \mcg{g}$ obtained by extending via the identity map 
along $D$. 
The kernel of this surjection is generated by two kinds of elements: the Dehn twist along the boundary curve, and ``push'' maps along elements of $\pi_1(S_{g,1})$ \cite{Birman1974}.  Note that any map in this kernel extends to $X_g$.  Hence, we are justified in defining the handlebody subgroup $\handle{g,1}$ of $\mcg{g,1}$ as the pullback of $\handle{g}$.  

The handlebody group arises naturally in a number of applications in 3-manifold topology, particularly through Heegaard splittings of 3-manifolds.  
Our second result  
in this paper is to compute $\rho_2(\handle{g,1})$, given in Theorem~\ref{theorem:morita image of handle}.

The authors would like to thank the referee for helpful comments and suggestions.
 
\section{The second Johnson-Morita map}
In this section we will describe Perron's formulation \cite{Perron2004} of the second Johnson-Morita representation.  We will give a precise characterization of the image of the mapping class group under this map.  First, it will be useful to review the image of the first Johnson-Morita representation, i.e., the symplectic group.  
\subsection{The symplectic group}
\label{SS:the symplectic group}
 
The group $H = H_1(S_{g,1})$ is free abelian with free basis $a_1,\dots, a_g,$ $b_1,\dots,b_g$, as in Figure~\ref{figure:curves}(a), and has a symplectic intersection form given by signed intersection of curves which is preserved by every mapping class $f\in \mcg{g,1}$.  In the basis above, the intersection form is given by the the matrix $J$ with $g \times g$ block form
\begin{equation}
\label{equation:J def}
J = \left ( \begin{array}{ccc} 0 &  -I \\ I & 0 \end{array} \right )
\end{equation}
The intersection form got by acting by the linear transformation $M$ on an intersection form with matrix $L$ is given by $ML\overline{M}$ where $\overline{M}$ denotes the transpose of $M$.  Hence for every $M$ in the image of the mapping class group
\begin{equation}
\label{equation:symplectic identity}
MJ\overline{M} = J, \ \ \ \mbox{ or equivalently } \ \ \ \overline{M}JM = J
\end{equation}
In fact (\ref{equation:symplectic identity}) is a sufficient condition for $M$ to be in the image of the mapping class group under $\rho_1$.  It is sometimes useful to write a symplectic matrix $M$ in $g \times g$ block form as
\begin{equation*}
M = \left( \begin{array}{ccc} S & T \\ P & Q \end{array} \right)
\end{equation*}

A convenient consequence of (\ref{equation:symplectic identity}) is that $M^{-1} = J\overline{M}J^{-1}$.  In block form this becomes
\begin{equation*}
\left( \begin{array}{ccc} S & T \\ P & Q \end{array} \right)^{-1} = \left( \begin{array}{ccc} \overline{Q} & -\overline{T} \\ -\overline{P} & \overline{S} \end{array} \right)
\end{equation*}

The group of such matrices form the \emph{symplectic group}.  Writing $M$ and $\overline{M}$ in $g \times g$ block form
\begin{equation*}
M= \left( \begin{array}{ccc} S & T \\ P & Q \end{array} \right) , \ \ \ \ \ \ \ \ 
\overline{M} = \left( \begin{array}{ccc} \overline{S} & \overline{P} \\ \overline{T} & \overline{Q} \end{array} \right)
\end{equation*}
we derive the {\it symplectic constraints}, which follow  directly from the condition in (\ref{equation:symplectic identity}):
\begin{eqnarray} 
\label{equation:symplectic conditions}
{\rm (i)} \ Q\overline{S}  -  P\overline{T} = I , \ \ 
{\rm (ii)} \  S\overline{T}  \  {\rm symmetric}, \ \ {\rm (iii)} \  P\overline{Q} \  {\rm symmetric}.
\end{eqnarray}
%%% In an earlier version you asked for a check that the preceding is correct. I checked it, and it is. 

\subsection{Perron's formulation of $\rho_2$}  
\label{section: invariants from the second Johnson-Morita homomorphism} 

%%%%%%%%%%%%%%%%%%%%%%%%%%%%%%%%%%%%%%%%%%%%%%%%%%%%%%%%%%%%%%%%%%%%%
% SUBSECTION Step I: The second Johnson-Morita representation
%%%%%%%%%%%%%%%%%%%%%%%%%%%%%%%%%%%%%%%%%%%%%%%%%%%%%%%%%%%%%%%%%%%%%
The {\em Torelli group} $\torelli{g,1}$ is the kernel of the symplectic representation $\rho_1 : \mathcal{M}_{g,1} \to \sp(H)$.
Johnson proved, in \cite{Johnson1980}, that the image of the Torelli group under $\rho_2$
is $\wedge^3 H$.  In the next section we will go a step further, and describe, in Theorem~\ref{theorem:morita image}, the image of the full mapping class group $\mathcal{M}_{g,1}$ under $\rho_2$ noting that Morita \cite[Theorem 4.8]{Morita1993a} has already identified this image as being finite index in $ \frac{1}{2} \wedge^3 H \rtimes \sp(H)$.
We begin by summarizing
%%% end of change
 Morita's explicit description of $\rho_2$ as given in \cite[Section 4]{Morita1993a}.  Consider the 2-step nilpotent group $$\Phi_2 = \left\{ (\eta, y) \bigg | \eta \in \frac{1}{2} \wedge^2 H , \mbox{  }y \in H \right\}$$ with multiplication in $\Phi_2$ given by $(\eta, y)(\nu, z) = (\eta + \nu + \frac{1}{2} y \wedge z, y+z )$.  It contains a subgroup of finite index which can be identified (see \cite[Sec. 5.5]{MKS}) with the second nilpotent quotient $\pi / \pi^{(2)} = \pi/[\pi,[\pi,\pi]]$ of our surface group via
the homomorphism $\phi_2 : \pi \to \Phi_2$
\begin{equation*}
\phi_2 ( \xi_i ) = (0,x_i)
\end{equation*}
where $\{ \xi_1, \cdots, \xi_{2g} \}$ generate $\pi = \pi_1(S_{g,1})$ and $\{ x_1,\cdots, x_{2g} \}$ is our basis for $H = H_1(S_{g,1})$ (see Figure~\ref{figure:curves}(a-b)).  The group $\Phi_2$ can be viewed as a subgroup of the Mal'cev completion of the nilpotent group $\pi/\pi^{(2)}$.  Any automorphism of $\pi / \pi^{(2)}$ extends to the Mal'cev completion and preserves $\Phi_2$ so we may think of $\mcg{g,1}$ as acting on $\Phi_2$ \cite[Proposition 2.5]{Morita1993a}.

In \cite[Section 3]{Morita1993a} Morita describes a function $\mcg{g,1} \to \hom (H,\frac{1}{2}\wedge^2 H)$.  An automorphism $f$ of $\Phi_2$ coming from an automorphism of the Mal'cev completion of $\pi / \pi^{(2)}$ can be specified by the images
\begin{equation*}
 f(0,x_i) = (w_i,h_i) \ \ \ \ \ \ w_i \in \frac{1}{2} \wedge^2 H, \ h_i \in H
\end{equation*}
for each $x_i$.  The homomorphism $\rho_1(f):H \to H$ given by $\rho_1(f)(x_i) = h_i$ is just the image of $f$ under the symplectic representation.  Johnson looks at the homomorphism $\tilde{\tau}_2(f):H \to \frac{1}{2} \wedge^2 H$ given by
\begin{equation*}
\tilde{\tau}_2(f)(x_i) = w_i
\end{equation*}
The function $\tilde{\tau}_2:\mcg{g,1} \to \hom (H,\frac{1}{2} \wedge^2 H)$ is a homomorphism when restricted to the kernel $\torelli{g,1}$ of the symplectic representation.  Johnson \cite[Theorem 1]{Johnson1980} identifies its image as $\wedge^3 H \subset \hom (H,\frac{1}{2} \wedge^2 H)$, where $x_i \wedge x_j \wedge x_k \in \wedge^3 H$ is understood to be the homomorphism
\begin{equation}\label{equation:wedge3 in hom}
(x_i \wedge x_j \wedge x_k) (y) = \langle y , x_k \rangle x_i \wedge x_j + \langle y , x_i \rangle x_j \wedge x_k + \langle y , x_j \rangle x_k \wedge x_i 
\end{equation}
where $\langle , \rangle $ gives the symplectic pairing for vectors in $H$.
The map $\torelli{g,1} \rightarrow \wedge^3 H \subset \hom (H,\frac{1}{2} \wedge^2 H)$ is usually referred to as the {\em Johnson homomorphism}.  

Morita \cite[Section 3]{Morita1993a} begins by considering this map $\tilde{\tau}_2:\mcg{g,1} \to \hom (H,\frac{1}{2} \wedge^2 H)$ (in Morita's notation this is the map $\tilde{k}$).   While not a homomorphism it is a crossed homomorphism with respect to the symplectic action of the mapping class group on $\hom (H,\frac{1}{2} \wedge^2 H)$.  In other words, the map $\tilde{\tau}_2$ satisfies:
\begin{equation*}
\tilde{\tau}_2(fg) = \tilde{\tau}_2(f) + \rho_1(f)\tilde{\tau}_2(g) \ \ \ \ \ f,g \in \mcg{g,1}
\end{equation*}
Choose $R \in \sp(H)$,  $y \in H$, and $m \in \hom (H,\frac{1}{2} \wedge^2 H)$.  We note that the action of $\sp(H)$ on $\hom (H,\frac{1}{2} \wedge^2 H)$ in the
equation above (and in the remainder of this paper) is the natural ``change-of-basis'' action:
\begin{equation}\label{equation:sp action on hom H wedge 2 H}
(Rm)(y) = R m (R^{-1} y) 
\end{equation}
The crossed homomorphism property is exactly what is needed for the map $\tilde{\rho}_2 : \mcg{g,1} \to \hom (H,\frac{1}{2} \wedge^2 H) \rtimes \sp(H)$
given by
\begin{equation*}
\tilde{\rho}_2(f) = (\tilde{\tau}_2 (f), \rho_1 (f))
\end{equation*}
to be a homomorphism. The homomorphism $\tilde{\rho}_2$ gives the action of $\mcg{g,1}$ on $\phi_2(\pi)\subset \Phi_2$, via
the action of $(r,R) \in \hom (H,\frac{1}{2} \wedge^2 H) \rtimes \sp(H)$ on $\Phi_2$:
\begin{equation}
\label{equation:action of morita image on pi image}
\left( r, R \right) * (\eta, y) =  \left( r(Ry)+ R\eta, Ry \right)
\end{equation}

Morita shows that by modifying the crossed homomorphism $\tilde{\tau}_2:\mcg{g,1} \to \hom (H,\frac{1}{2} \wedge^2 H)$, one obtains a crossed homomorphism $\tilde{\tau}_2'$ (Morita denotes this map by $\tilde{k}'$ in \cite[Section 4]{Morita1993a} and $\tilde{k}$ in \cite[Section 5]{Morita1993a}) from $\mcg{g,1}$ to the submodule $\frac{1}{2}\wedge^3 H$ of $\hom (H,\frac{1}{2} \wedge^2 H)$ which extends the Johnson homomorphism.  We will modify $\tilde{\tau}_2$ to get a different crossed homomorphism $\tau_2:\mcg{g,1} \to \frac{1}{2}\wedge^3 H$ extending the Johnson homomorphism.  Our map $\tau_2$ is a trivial modification of Morita's map $\tilde{\tau}_2'$ which will lend itself to later calculations.

For any $m \in \hom (H,\frac{1}{2} \wedge^2 H)$, the map $\sigma_m: \mcg{g,1} \to \hom (H,\frac{1}{2} \wedge^2 H)$ given by
\begin{equation*}
\sigma_m(f) = m - \rho_1(f)m
\end{equation*}
is a crossed homomorphism.  Such a crossed homomorphism is called {\em principal}; two crossed homomorphisms are cohomologous if they differ by a principal crossed homomorphism \cite[Chapter IV.2]{Brown1994}.

Let $\kappa \in \hom (H,\frac{1}{2} \wedge^2 H)$ 
be the homomorphism
\begin{equation*}
\kappa (a_i)  =   \frac{1}{2}a_i \wedge b_i \ \ \ \ \ \ \kappa (b_i)  =   - \frac{1}{2}a_i \wedge b_i 
\end{equation*}
or equivalently
\begin{equation}\label{equation:kappa}
\kappa (x_i)  =   \frac{1}{2}x_i \wedge Cx_i 
\end{equation}
where $C$ is the $2g \times 2g$ matrix with $g \times g$ block form 
$\left( \begin{array}{cc} 0 & I \\ I & 0 \end{array} \right)$.
Define
\begin{equation}
\label{equation:good crossed hom}
\tau_2(f) = \tilde{\tau}_2(f) + \kappa - \rho_1(f)\kappa
\end{equation} 

This is the crossed homomorphism that Perron \cite[Remark 5.5]{Perron2004} denotes $-\frac{1}{6}\widetilde{A_1}$.
We note that by comparing the above with \cite[Proposition 4.7]{Morita1993a}, it is straightforward to see that Morita's crossed homomorphism 
$\tilde{\tau}'_2$ can be expressed as
\begin{equation*}
\tilde{\tau}'_2(f) = \tau_2(f) + m - \rho_1(f)m
\end{equation*}
where $m = - \frac{1}{2}(\sum_{i=1}^g a_i + b_i) \wedge (\sum_{i=1}^g a_i \wedge b_i)$.  In other words, the map $\tau_2$ and Morita's original map $\tilde{\tau}'_2$ are cohomologous, that is, they represent the same element of $H^1(\mcg{g,1},\frac{1}{2}\wedge^3H)$.

We can now define a homomorphism $\rho_2:\mcg{g,1} \to \frac{1}{2}\wedge^3 H \rtimes \sp(H)$ as follows:
\begin{equation*}
\rho_2(f) = (\tau_2 (f), \rho_1 (f))
\end{equation*}
Using (\ref{equation:good crossed hom}), (\ref{equation:action of morita image on pi image}), (\ref{equation:sp action on hom H wedge 2 H}), and (\ref{equation:wedge3 in hom}), we obtain the correct action of $\rho_2(\mcg{g,1})$ on $\Phi_2$:
\begin{eqnarray}
\textstyle
\lefteqn{ \left( \sum r_{ijk} x_i \wedge x_j \wedge x_k  ,  R \right) * (\eta, y)} \nonumber \\
&  =  &
\label{equation:easy action}
\left( R \eta - \kappa(Ry) + R(\kappa(y)) + r(y),  Ry \right) \\ 
& = & \scriptstyle \left( R \eta - \kappa(Ry) + R(\kappa(y))  + \sum r_{ijk} \left( \begin{array}{c} \scriptstyle \langle Ry , x_k \rangle x_i \wedge x_j \\
\scriptstyle + \langle Ry, x_i \rangle x_j \wedge x_k \\ \scriptstyle + \langle Ry , x_j \rangle x_k \wedge x_i \end{array} \right),  Ry \right)
\label{equation:action of perron image on pi image}
\end{eqnarray}
where $\langle , \rangle $ is the symplectic pairing on $H$ and the sums are taken over $1 \leq i < j < k \leq 2g$.

\subsection{Calculating the image of the mapping class group}
In this section we compute $\rho_2(\mathcal{M}_{g,1})$.   See Theorem~\ref{theorem:morita image} below.

Recall the map $\phi_2: \pi \rightarrow \Phi_2$ given in the previous section.  It will be helpful for us to identify  $\phi_2(\pi) \subset \Phi_2$ precisely.  The gist of the following lemma is that for pairs in the image of $\phi_2$, the second coordinate determines the first coordinate modulo 1.  
\begin{lemma}
\label{lemma:pi image 2}
The image of $\pi$ under the map $\phi_2$ is given as follows.  
$$\phi_2(\pi) = \left\{ \left( \sum_{1 < i < j < 2g} \left( n_{ij} + 
\frac{l_il_j}{2} \right) x_i \wedge x_j \, , \,\, \sum_{i=1}^{2g} 
l_i x_i \right ) \Bigg| \mbox{  } n_{ij},\mbox{ }l_i \in \ints \right\}$$ 
\end{lemma}

\begin{proof}
Let $G \subset \Phi_2$ denote the set on the right-hand side of the equation in the lemma.
We claim that the set
$G$ is a subgroup of $\Phi_2$. 
First, $G$ is closed under inversion 
since $(\eta, y)^{-1} = (-\eta,-y)$.  For closure under products 
consider
\begin{eqnarray*}
\lefteqn{ \left( \sum_{1 < i < j < 2g} \left( n_{ij} + 
\frac{l_il_j}{2} \right) x_i \wedge x_j \, , \,\, \sum_{i=1}^{2g} 
l_i x_i \right ) } \\
& & \cdot \left( \sum_{1 < i < j < 2g} \left( n_{ij}' + 
\frac{l_i'l_j'}{2} \right) x_i \wedge x_j \, , \,\, \sum_{i=1}^{2g} 
l_i' x_i \right ) \\
& = & \left( \sum_{1 < i < j < 2g} \left( \begin{array}{c} \scriptstyle n_{ij} + n'_{ij} + 
\frac{l_il_j}{2} \\ \scriptstyle+ \frac{l_i'l_j'}{2} + \frac{l_il'_j}{2} 
- \frac{l_jl_i'}{2} \end{array} \right) x_i \wedge x_j \, , \,\, 
\sum_{i=1}^{2g} (l_i + l'_i)x_i \right )
\end{eqnarray*}
This product is in $G$ because $l_il_j + l_i'l_j' + l_il'_j - l_jl_i' 
\equiv (l_i + l'_i)(l_j + l'_j) \mod 2$. 

Clearly, $G$ contains each generator $\phi_2(\xi_i) = (0,x_i)$ of $\phi_2(\pi)$. 
For the reverse inclusion, note that any element of the form
\begin{equation*}(0,x_i)(0,x_j)(0,-x_i)(0,-x_j) 
= (x_i \wedge x_j, 0)
\end{equation*} 
lies in $\phi_2(\pi)$.  In fact such an element is in the 
center of $G$.  Now, any element of $G$ can be written as a 
product of $(0,x_i)$'s to get the correct second coordinate, followed 
by a product of $(x_i \wedge x_j, 0)$'s to get the correct first coordinate. Hence 
$G \subset \phi_2(\pi)$.
\end{proof}

We are almost ready to 
characterize the subgroup $\rho_2(\mcg{g,1}) \subset  \frac{1}{2} \wedge^3 H \semi \sp(H)$. We begin with a simple yet fundamental observation.

\begin{remark}\label{remark:differing by integers}
Suppose $R$ is a symplectic matrix and $(r_1,R), 
(r_2,R) \in \rho_2(\mcg{g,1})$.  Then $(r_1,R)^{-1} = 
(-R^{-1}r_1,R^{-1}) \in \rho_2(\mcg{g,1})$ so
\begin{equation*}
(r_2,R)(-R^{-1}r_1,R^{-1}) = (r_2-r_1,I) \in \rho_2(\mcg{g,1}).
\end{equation*}
In other words, we have that $(r_2-r_1,I) \in \rho_2(\torelli{g,1})$.  Using  
Johnson's characterization of $\tau_2(\torelli{g,1})$ 
\cite[Theorem 1]{Johnson1980} we conclude that if two elements of 
$\rho_2(\mcg{g,1})$ have identical symplectic matrices, then their 
$\frac{1}{2} \wedge^3 H$ coordinate must differ by an 
\emph{integral} element of $\wedge^3 H$.  
\end{remark}

As a consequence of this observation, we expect 
that the symplectic matrix $R$ will determine the coefficients of $r_1$ 
and $r_2$ modulo $1$.  Theorem~\ref{theorem:morita image} makes this 
precise and gives the characterization of $\rho_2(\mcg{g,1})$.  First we give a short definition.
\begin{definition}
Given three $n$-dimensional vectors $\vec{w}=(w_1,\dots, w_n)$, 
$\vec{y}=(y_1,\dots, y_n)$, $\vec{z}=(z_1,\dots, z_n)$ in basis
$\mathcal{B}$,
their \emph{$\mathcal{B}$-triple dot product} is the scalar 
$$\bullet_\mathcal{B}(\vec{w},\vec{y},\vec{z}) = 
\sum_{i=1}^n w_i y_i z_i.$$ When the basis $\mathcal{B}$ is clear, we 
will write $\bullet (\vec{w},\vec{y},\vec{z})$.
\end{definition}

Recall that $J$ is the matrix given in (\ref{equation:J def}). 
\begin{theorem}
\label{theorem:morita image}
Let $R \in \sp (2g, \ints)$ be an arbitrary symplectic matrix.  Let 
$r$ be any element of $\frac{1}{2} \wedge^3 H$ with 
$r = \sum_{1\leq i < j < k \leq 2g} r_{ijk} x_i \wedge x_j \wedge x_k$. 
Then  $(r,R) \in \rho_2(\mcg{g,1})$ if and only if
$$r_{ijk} \equiv \frac{E_{ijk}}{2} \mod 1$$
where
\begin{eqnarray*}
E_{ijk} & = & \bullet(\row_i(RJ), \row_j(R),\row_k(R)) \\
& & - \bullet(\row_i(R), \row_j(RJ),\row_k(R)) \\
& & + \bullet(\row_i(R), \row_j(R),\row_k(RJ))
\end{eqnarray*}
for all $1\leq i < j < k \leq 2g$. 
\end{theorem}

\begin{proof}
Let $(r,R) \in \rho_2(\mcg{g,1})$, and let
\begin{equation*}
r = \sum_{1\leq i < j < k \leq 2g} r_{ijk} x_i \wedge x_j \wedge x_k.
\end{equation*}
For $1 \leq i,j,k \leq 2g$ we set $r_{ijk} = 0$ unless $i < j < k$.
The group $\rho_2(\mcg{g,1})$ preserves 
$\phi_2(\pi)$, described in Lemma~\ref{lemma:pi image 2}.  
Let $x_n$ be an arbitrary basis element of $H$, and consider the action of $(r, R)$ on $(0, x_n)$.  We will use the standard notation $M_{ij}$ to denote the entry in the $i^{th}$ row and $j^{th}$ column of a matrix $M$ throughout.  
By (\ref{equation:action of perron image on pi image}), we get that the second coordinate of $(r,R)*(0,x_n)$ is simply $Rx_n$, which we can write as  $\sum_{i=1}^{2g} R_{in} x_i$, with an eye on eventually applying Lemma~\ref{lemma:pi image 2}.  
Using (\ref{equation:action of perron image on pi image}) and (\ref{equation:kappa}), we obtain the following for the first coordinate of $(r,R)*(0,x_n)$:
\begin{equation*} - \kappa(R x_n) + R(\kappa (x_n)) + \sum_{1\leq i < j < k \leq 2g} r_{ijk} \left( \begin{array}{c}
  \langle Rx_n , x_k \rangle x_i\wedge x_j \\ + \langle Rx_n , x_i \rangle x_j\wedge x_k \\
  - \langle R x_n, x_j \rangle x_i\wedge x_k \end{array} \right) 
\end{equation*}
Notice that under the symplectic pairing $\langle Rx_n , x_k \rangle = (JR)_{kn}$ so
the above can be rewritten:
\begin{eqnarray*}
&   &  - \kappa \left(\sum_{i=1}^{2g} R_{in} x_i \right) + R\left(\frac{1}{2} x_n\wedge Cx_n \right) \\
&   &  + \sum_{1\leq i < j < k \leq 2g} r_{ijk} \left( \begin{array}{c}
  ((JR)_{kn}) x_i\wedge x_j \\ + ((JR)_{in}) x_j\wedge x_k \\
  - ((JR)_{jn}) x_i\wedge x_k \end{array} \right) \\
& = & - \left(\sum_{i=1}^{2g} \frac{R_{in}}{2} x_i \wedge Cx_i \right) + \left(\sum_{1 \leq i,j \leq 2g} 
  \frac{R_{in} (RC)_{jn}}{2} x_i \wedge x_j \right) \\
&   &  + \sum_{1\leq i < j < k \leq 2g} r_{ijk} \left( \begin{array}{c}
  ((JR)_{kn}) x_i\wedge x_j \\ + ((JR)_{in}) x_j\wedge x_k \\
  - ((JR)_{jn}) x_i\wedge x_k \end{array} \right) \\
& = & \left(\sum_{i=1}^{g} \frac{(CR)_{in} - R_{in}}{2} x_i \wedge x_{i+g} \right) \\
&   & + \left(\sum_{1 \leq i < j \leq 2g}  
  \frac{R_{in} (RC)_{jn} - R_{jn} (RC)_{in}}{2} x_i \wedge x_j \right)  \\
&   &  + \sum_{1\leq i < j < k \leq 2g} r_{ijk} \left( \begin{array}{c}
  ((JR)_{kn}) x_i\wedge x_j \\ + ((JR)_{in}) x_j\wedge x_k \\
  - ((JR)_{jn}) x_i\wedge x_k \end{array} \right)
\end{eqnarray*}
Now, applying Lemma 
\ref{lemma:pi image 2} to the coefficient of $x_p\wedge x_q$, where 
$p < q$, gives
\begin{eqnarray*}
\lefteqn{\frac{\delta_{q,p+g}((CR)_{pn} - R_{pn}) + R_{pn} (RC)_{qn} - R_{qn} (RC)_{pn}}{2} } \\
& & + \sum_{i=1}^{2g} \left( r_{ipq} (JR)_{in} - r_{piq} (JR)_{in} 
  + r_{pqi} (JR)_{in} \right) \equiv \frac{R_{pn}R_{qn}}{2} \mod 1
\end{eqnarray*}
Note that for fixed $i,p,q$, at most one of the $r$-coefficients in the above summation is nonzero.  For bookkeeping purposes, when $1 \leq j<r\leq 2g$ we define $\vec{r}_{jk}$ be the $2g$-dimensional column vector whose $i^{th}$ 
entry is $r_{ijk}$ if $i < j$, $-r_{jik}$ if $j < i < k$, 
$r_{jki}$ if $k < i$, and $0$ otherwise. If $\mbox{col}_n(M)$ denotes the $n^{th}$ column vector of $M$, we may rewrite this to obtain that $\mbox{col}_n(JR) \cdot \vec{r}_{pq}$ is congruent $(\mod 1)$ to
\begin{equation*}
\frac{\delta_{q,p+g}(R_{pn} - (CR)_{pn}) + R_{pn}R_{qn} - R_{pn} (RC)_{qn} + R_{qn} (RC)_{pn}}{2}
\end{equation*}
In order to write this a bit more compactly, for $1 \leq j < k \leq 2g$, we define $\vec{t}_{jk}$ to be the $2g$-dimensional 
column vector whose $i^{th}$ entry is $\delta_{k,j+g}(R_{ji}- (CR)_{ji}) + R_{ji} R_{ki} - R_{ji} (RC)_{ki} + R_{ki} (RC)_{ji}$.
Combining the equations above for all $1\leq n \leq 2g$ we get:
\begin{equation*}
\overline{JR} \vec{r}_{pq}  \equiv  \frac{\vec{t}_{pq}}{2} \mod 1 \ \ \ \ \ \forall 1 \leq p < q \leq 2g
\end{equation*}
Solving for $\vec{r}_{pq}$, we obtain:
\begin{equation*}
\vec{r}_{pq}  \equiv \frac{(\overline{JR})^{-1} \vec{t}_{pq}}{2} 
\mod 1 
\end{equation*}
Since $R$ is assumed to be symplectic, we can rewrite this as: 
\begin{equation*}
\vec{r}_{pq}  \equiv \frac{RJ\vec{t}_{pq}}{2} \mod 1
\end{equation*}
Observe that the $i^{th}$ entry of the vector on the right-hand side is 
\begin{eqnarray}
& & \frac{1}{2}\delta_{q,p+g} \row_i(RJ) \cdot (\row_p(R) - \row_p(CR)) \nonumber \\
& & +\frac{1}{2} \bullet(\row_i(RJ), \row_p(R),\row_q(R)) \nonumber \\
& & -\frac{1}{2} \bullet(\row_i(RJ), \row_p(R),\row_q(RC)) \nonumber \\
\label{equation:with delta}
& & +\frac{1}{2} \bullet(\row_i(RJ), \row_p(RC),\row_q(R))
\end{eqnarray}
We are interested in calculating the coefficients $r_{ipq}$ for $1 \leq i < p < q \leq 2g$.  Thus we are interested
in the $i^{th}$ entry of $\vec{r}_{pq}$ when $1 \leq i < p < q \leq 2g$.  If $q\neq p+g$ then $\delta_{q,p+g} =0$.
Assume that $q=p+g$. Then $1 \leq i < p \leq g$, and if we write $R = \left( \begin{array}{ccc} S & T \\ P & Q \end{array} \right)$, we have
\begin{eqnarray*}
\lefteqn{\row_i(RJ) \cdot (\row_p(R) - \row_p(CR)) } \\
& = & \row_i(T) \cdot \row_p(S) - \row_i(S) \cdot \row_p(T) \\
&   & - \row_i(T) \cdot \row_p(P) + \row_i(S) \cdot \row_p(Q) \\
& = & (T\overline{S})_{ip} - (S\overline{T})_{ip} - (T\overline{P})_{ip} + (S\overline{Q})_{ip} \\
& = & 0 - 0
\end{eqnarray*}
The last equality results from using the symplectic conditions (\ref{equation:symplectic conditions} i,ii) and by our assumption that $i \neq p$.  Thus we may drop the first term of (\ref{equation:with delta}).  In other words, for $1 \leq i < p < q \leq 2g$ the $i^{th}$ entry of $\vec{r}_{pq}$ (mod 1) is given by
\begin{eqnarray*}
& & \frac{1}{2} \bullet(\row_i(RJ), \row_p(R),\row_q(R)) \\
& & -\frac{1}{2} \bullet(\row_i(RJ), \row_p(R),\row_q(RC)) \\
& & +\frac{1}{2} \bullet(\row_i(RJ), \row_p(RC),\row_q(R)) \ \mod 1
\end{eqnarray*}
For aesthetic reasons we rewrite the expression above more symmetrically
to show that $i^{th}$ entry of $\vec{r}_{pq}$ (mod 1) is:
\begin{eqnarray*}
& & \frac{1}{2} \bullet(\row_i(RJ), \row_p(R),\row_q(R)) \\
& & -\frac{1}{2} \bullet(\row_i(R), \row_p(RJ),\row_q(R)) \\
& & +\frac{1}{2} \bullet(\row_i(R), \row_p(R),\row_q(RJ)) \ \mod 1
\end{eqnarray*}
We have just shown 
that the $\binom{2g}{3}$ equations in the statement of the lemma are 
necessary for $(r,R)$ to be an element of $\rho_2(\mcg{g,1})$.  
Since the symplectic representation $\rho_1$ is surjective,
$\rho_2(\mcg{g,1})$ contains an element of the form $(r, R)$ for any given $R$.  Johnson \cite[Theorem 1]{Johnson1980} showed that any element of the 
form $(w,I)$ with $w \in \wedge^3 H$ is in $\rho_2(\mcg{g,1})$.  Then
if $(r,R) \in \rho_2(\mcg{g,1})$, so is $(w,I)(r,R) = (w + r, R)$ for any $w \in \wedge^3 H$. 
Hence we can hit any other possible choice of the coefficients $r_{ijk}$ satisfying the ``mod 1'' conditions imposed by $R$ by composing our map with different choices of Torelli elements.  This shows sufficiency. 
\end{proof}

%%%%%%%%%%%%%%%%%%%%%%%%%%%%%%%%%%%%%%%%%%%%%%%%%%%%%%%%%%%%%%%%%%%%%
% SUBSECTION Step II: The symplectic image of the handlebody group
%%%%%%%%%%%%%%%%%%%%%%%%%%%%%%%%%%%%%%%%%%%%%%%%%%%%%%%%%%%%%%%%%%%%%

\section{The handlebody group}
Our primary goal in this section is to compute $\rho_2( \handle{g,1} )$ explicitly.  We will begin with some known algebraic characterizations of $\handle{g,1}$ and of $\rho_1(\handle{g,1})$ which will be helpful to us, and use them to derive an analogous characterization at the second level.  Thus equipped, we derive an explicit formulation of $\rho_2(\handle{g,1})$ in Section~\ref{subsection:image of handlebody}.  

\subsection{Algebraic characterizations of the handlebody subgroup}
\label{SS:symplectic image of handlebody}

Let $\mathfrak{b}$ denote the normal 
closure in $\pi$ of $\{ \beta_1,\dots,\beta_g \}$.  Note that $\mathfrak{b}$ is also the kernel 
of the homomorphism $\pi \to \pi_1 (X_g)$ induced by inclusion.

The following proposition was first proved by McMillan \cite{McMillan}. 
The proof given here was suggested to the authors by Saul Schleimer. 
\begin{proposition} 
\label{proposition:characterization of handlebody}
The handlebody subgroup $\handle{g,1}$ of the mapping class group 
$\mcg{g,1} \subset \aut (\pi_1(S_{g,1}))$ is precisely the subgroup which 
preserves $\mathfrak{b}$. 
\end{proposition}

\begin{proof}
One direction is immediate; in order for a mapping class in 
$\mcg{g,1}$ to extend to the $X_g$ it must preserve $\mathfrak{b}$.  
Now suppose $f$ is a mapping class which preserves $\mathfrak{b}$.
Then $f$ sends each $\beta_i$ to a loop that can be represented
by a simple closed curve which is trivial in $\pi_1(X_g)$.  
Dehn's Lemma \cite{Papa1957} shows that these curves bound disks in 
$X_g$ that can be made disjoint.  By matching these disks to
the ones bounded by each $\beta_i$ we may 
construct a homeomorphism from $X_g$ to itself restricting to $f$
on its boundary.
\end{proof}

Moving on to level one of the Johnson-Morita representations, Birman has shown that the image of the handlebody group in 
$\sp(2g, \ints)$ is particularly nice \cite[Lemma 2.2]{Birman1975}.  All subblocks are $g \times g$ matrices.
\begin{proposition}[Birman]
\label{proposition:sphandlebody image}
The image of the handlebody group under the symplectic representation is characterized by a $g \times g$ block of zeroes in the upper-right corner.  That is, $$\rho_1(\handle{g,1}) = \left\{ M \in \sp(2g;\ints)
\bigg| M \mbox{ has block form } \left( \begin{array}{ccc} * & 0 
\\ * & * \end{array} \right) \right\}$$
\end{proposition}

Sufficiency is shown 
in \cite{Birman1975} by exhibiting generators for $\rho_1(\handle{g,1})$
which are in the image of the handlebody group.
The necessity of this condition for membership in  
$\rho_1(\handle{g,1})$ follows from the observation that 
in the handlebody $X_g$, the homology classes of the generators of
type $b_i$ are all 0.  Any homeomorphism of $S_g$ which
extends to $X_g$ must take trivial
elements in the homology of the handlebody to trivial elements 
in the homology of the handlebody.  In other words, $\rho_1(\handle{g,1})$ is characterized by the property that its elements must preserve the subgroup of $H$
generated by the $b_i$'s. 

We will now give a second-level analogue of these characterizations by describing a subgroup of $\pi / \pi^{(2)}$ which must be preserved by $\rho_2 (\handle{g,1})$, thus giving a restriction on the image of 
the handlebody group.

The second Johnson-Morita homomorphism is given by the action of 
$\mcg{g,1}$ on the nilpotent quotient $\pi / \pi^{(2)}$.
Let $\mathfrak{b} \subset \pi$ be as above, and recall from Section~\ref{section: invariants from the second Johnson-Morita homomorphism} the map $\phi_2: \pi \rightarrow \Phi_2$ be as above. The 
following lemma computes $\phi_2 (\mathfrak{b})$.

\begin{lemma}
\label{lemma:phi2 of pi} 
\begin{equation*}
\phi_2 (\mathfrak{b}) = \left\{ \left( \begin{array}{c} \scriptstyle \sum_{1 \leq i,j \leq g} 
m_{ij} a_i \wedge b_j \\ \scriptstyle + \sum_{1\leq i < j \leq g} \left( 
n_{ij} + \frac{l_i l_j}{2}\right)b_i \wedge b_j \end{array} \, , \,\,
\sum_{i=1}^{g} l_i b_i \right) \bigg| m_{ij},n_{ij}, l_i \in \ints 
\right\}
\end{equation*}
\end{lemma}
\begin{proof}
In light of Lemma~\ref{lemma:pi image 2}, the right-hand side above is clearly the 
kernel of the quotient homomorphism $\pi / \pi^{(2)} \to 
\pi_1(X_g) / \pi_1(X_g)^{(2)}$.
\end{proof}

Now that we have identified $\phi_2(\mathfrak{b})$ we will describe 
$\rho_2(\handle{g,1})$.  

\subsection{Image of the handlebody subgroup under $\rho_2$}
\label{subsection:image of handlebody}

Theorem~\ref{theorem:morita image} above gives  $\rho_2(\mcg{g,1})$.   The missing ingredient for a characterization of $\rho_2(\handle{g,1})$ is $\rho_2(\torelli{g,1} \cap \handle{g,1})$ which was computed by Morita.
\begin{proposition}[{\cite[Lemma~2.5]{Morita1989}}] 
\label{proposition:tau2 of handlebody intersect torelli} 
$\rho_2(\torelli{g,1} \cap \handle{g,1})$ is the free abelian group 
with free basis:
$$ (b_i \wedge b_j \wedge b_k,I),   \ \ (a_i\wedge b_j\wedge b_k,I)
, \ \ {\rm and} \ \ (a_i\wedge a_j \wedge b_k,I) \ \  
1\leq i,j,k \leq g. $$
\end{proposition} 

Now we have the tools to assemble a description of  
$\rho_2( \handle{g,1} )$.  The following theorem gives a complete characterization of $\rho_2(\handle{g,1})$; it says that an element is in this image if and only if its first factor has no ``triple-a'' terms and its second factor has the form of Proposition~\ref{proposition:sphandlebody image}.  

\begin{theorem}
\label{theorem:morita image of handle}
Let $R \in \sp (2g, \ints)$ be an arbitrary symplectic matrix.  Let 
$r$ be any element of $\frac{1}{2} \wedge^3 H$ with 
$r = \sum_{1\leq i < j < k \leq 2g} r_{ijk} x_i \wedge x_j \wedge x_k$. 
Then  $(r,R) \in \rho_2(\handle{g,1})$ if and only if all of the 
following three conditions hold:
\begin{enumerate}
\item $R$ has $g \times g$ block form
$\left( \begin{array}{cc} * & 0 \\ * & * \end{array} \right)$
\item
$r_{ijk} \equiv \frac{1}{2}E_{ijk} \mod 1$
for all $1\leq i < j < k \leq 2g$. 
\item $r_{ijk} = 0$ for all $i,j,k$ with $0 \leq i < j < k \leq g$. 
(i.e. $r$ contains no terms of the form $a_i \wedge a_j \wedge a_k$.)
\end{enumerate}
\end{theorem}
We refer the reader to Theorem~\ref{theorem:morita image} for
the definition of $E_{ijk}$, which depends on the matrix $R$.
\begin{proof}
The necessity of condition 1 has already been established in \cite[Lemma 2.2]{Birman1975}.
We claim that only elements of $\frac{1}{2}\wedge^3 H 
\rtimes \sp (H)$ satisfying condition 3 above preserve $\phi_2(\mathfrak{b})$ under
the action of (\ref{equation:action of perron image on pi image}).
Suppose $R$ is symplectic with the required block form and $r$ 
contains a term of the form $c a_i \wedge a_j \wedge a_k$.  Since $R^{-1}$ must satisfy condition 1 above and using Lemma~\ref{lemma:phi2 of pi}, 
there is an element $(\nu,R^{-1}b_i) \in \phi_2(\mathfrak{b})$ where
$\nu$ has only terms of the form $\frac{1}{2}b_n \wedge b_m$. Applying  
(\ref{equation:easy action}) we get
\begin{eqnarray*}
\lefteqn{(r,R)*(\nu,R^{-1}b_i) = } \\
& = & \left( R(\nu) + \kappa ( R R^{-1}b_i) + R \kappa ( R^{-1}b_i ) + r (R R^{-1}b_i), R R^{-1} b_i \right) \\
& = & \left( R(\nu) + \kappa ( b_i) + R \kappa ( R^{-1}b_i ) + r (b_i), b_i \right)
\end{eqnarray*}
Consider each of the terms in the first coordinate of the ordered pair above.
Since $\nu$ only has terms of the form $\frac{1}{2}b_n \wedge b_m$ and
the matrix $R$ has the block form given in condition 1, we must have that $R(\nu)$
contains no terms of the form $a_j \wedge a_k$.  The image of the homomorphism
$\kappa$ has no $a_j \wedge a_k$ terms so neither $\kappa ( b_i)$ nor 
$\kappa ( R^{-1}b_i )$ contains any $a_j \wedge a_k$ terms.  Application 
of the matrix $R$ preserves this quality; hence $R \kappa ( R^{-1}b_i )$
contains no $a_j \wedge a_k$ terms.  We can see using (\ref{equation:wedge3 in hom}) that $r (b_i)$ will
contain a term of the 
form $-c a_j \wedge a_k$ by construction.  Then Lemma~\ref{lemma:phi2 of pi} implies that $c=0$.  It follows that the two conditions of the corollary are necessary.
 
For each $R$ satisfying 1 there is some mapping class $f \in 
\handle{g,1}$ with $\rho_1(f) = R$ as shown in \cite[Lemma 2.2]{Birman1975}.  
We have shown that $\rho_2(f)$ satisfies conditions 1 and 2.  Applying Proposition~\ref{proposition:tau2 of handlebody intersect torelli} we can get 
every other element of the form $(w,R)$ satisfying 1 and 2 as a 
product $(z,I)\rho_2(f)$ where $(z,I) \in \rho_2(\torelli{g,1} \cap 
\handle{g,1})$.  This establishes sufficiency. 
%Finally, the last statement of the theorem follows if one notes that 
%if the upper right $g \times g$ block of $R$ is $0$ then $r_{{}_R}$ 
%has no terms of the form $a_i \wedge a_j \wedge a_k$.
\end{proof}

%%%%%%%%%%%%%%%%%%%%%%%%%%%%%%%%%%%%%%%%%%%%%%%%%%%%%%%%%%%%%%%%%%%%%
%%%%%%%%%%%%%%%%%%%%%%%%%%%%%%%%%%%%%%%%%%%%%%%%%%%%%%%%%%%%%%%%%%%%%
% BIBLIOGRAPHY
%%%%%%%%%%%%%%%%%%%%%%%%%%%%%%%%%%%%%%%%%%%%%%%%%%%%%%%%%%%%%%%%%%%%%
%%%%%%%%%%%%%%%%%%%%%%%%%%%%%%%%%%%%%%%%%%%%%%%%%%%%%%%%%%%%%%%%%%%%%

\noindent \small \em
{\rm Joan S. Birman}\\
Department of Mathematics\\
Barnard College of Columbia University \\
New York, New York 10027\\
USA\\
\\
{\rm Tara E. Brendle}\\
Department of Mathematics\\
Louisiana State University\\
Lockett Hall\\
Baton Rouge, LA  70803\\
USA\\
\\
{\rm Nathan Broaddus}\\
Department of Mathematics\\
University of Chicago\\
5734 S. University Ave.\\
Chicago, IL 60637\\
USA

\end{document}